# Transforming Mathematics

Benjamin Braun

Professor of Mathematics

University of Kentucky

benjamin.braun@uky.edu

https://sites.google.com/view/braunmath/

I believe that the development of AI tools will lead to a transformation of mathematics, both professionally and culturally. However, it will be a tragic missed opportunity if the massive shifts in culture and practice that are underway are motivated and influenced solely by AI. We must aspire to more than this. We can and should collaboratively strive for a future where our cultures, practices, and values are grounded first and foremost in empathy, compassion, justice, and humanity.

In this article, I call on every member of the mathematical community, whether you view yourself as a researcher, teacher, scholar, student, professional, hobbyist, or a combination of these and more, to take the following actions.

**Do not rush.** There is tremendous pressure on everyone involved in mathematics at this moment, pressure due to the increasing capabilities of AI models. Do not yield to this pressure. We can make reasonable temporary adjustments to handle short-term disruptions, and we do not have to make fast, ill-considered long-term decisions. The large AI companies have a financial incentive to convince us to make decisions that utilize their systems, and to do so in ways that risk substantial long-term harm to mathematical learning, knowledge, and communities. Professional societies, department leaders, parents, students, teachers, and faculty should collectively demand that we make thoughtful, careful, deliberate decisions.

**Expand your connections and listen to others**. We must collectively transcend the boundaries we have set for ourselves and seek a greater perspective and vision, across all our interests and activities in the mathematical sciences. Together, we must expand our view beyond our own priorities, our own institutions, and our own countries and continents. We cannot afford to hide our conversations within our silos of comfort, because decisions within each part of the world of mathematics will have ripple effects that extend far beyond our own interests and institutions. The decisions we make about how to handle AI in our research communities must be aligned with the decisions we make about how AI is handled in the education of students from preschool through graduate school, as all of these are interdependent.  If you are a department chair, if you are a committee member for a professional society, if you are in any position of leadership, create a working group that facilitates dialogue across the mathematical sciences. If you are not in a position of leadership, become a leader by creating opportunities for dialogue, reflection, and debate. We must each lead, and we must each follow.

**Learn about the practices of mathematics.** Invest time to pause and reflect: what are our collective mathematical practices? There exists an entire discipline devoted to this question, the philosophy of mathematical practice, which includes researchers in philosophy, mathematics, mathematics education, psychology, and other disciplines. There are a wide variety of practices that are important to mathematical thought and activity, many of which we never teach explicitly, whether in classes or through informal mentoring and guidance, and some of which we are often not consciously aware of. Ask yourself: when you are collaborating on a mathematical task, do you make intentional decisions about how you collaborate? Or do you simply collaborate, without an explicit knowledge of your collaborative actions? When you learn something new, or teach something to others, are you consciously aware of how your actions and behaviors impact learning and feelings of belonging? Think about all the things that we do when we are “doing mathematics”, from reading to writing to experimenting to mentoring to speaking and beyond. How will our decisions about AI influence each of those practices, each of those habits, whether conscious or subconscious? This is too much to unpack alone, and it requires time to understand. Before we can make permanent decisions about AI, we must deeply reflect on what our decisions will do, both for us and to us, and we must engage in collaborative reflection.

**Embrace our obligation to support others**. Now is not the time to prioritize prestige or self-centered interests. Our responsibility is to think beyond ourselves, to those throughout our mathematical communities, and to turn our thoughts toward the needs of those who are most marginalized within the worlds of mathematics. Our focus must be foremost on protecting the next generation of mathematical learners, ensuring that they can develop their unique mathematical identities and share power in the mathematical community. It is our collective responsibility to support young people, and to welcome them into a culture of humane and human mathematics. This includes students at every level, from graduate study to preschool and primary school. If you are a teacher or a professor, find a student and talk with them about mathematics, and encourage them to engage in the struggle to learn something new.

**Insist that AI tools be ethically developed and implemented.** If AI tools become part of our mathematical practices, we must ensure that we are using them in an ethical manner and that we use tools that are developed and implemented ethically. Many of the current practices of AI companies are neither ethical nor just. We should ensure that our mathematical work and the tools that we use for it do not rely on the oppression or exploitation of other people or the natural world. Before you commit to using AI tools in your practices of mathematics, learn about the tools. What data were they trained on? How has the company that produced the tool treated its workers? How has it treated the communities in which it operates? Activate your inner skeptic and view these tools and their makers with a critical lens. We do not have to use tools that have caused or will cause harm to others.

**Do not believe in myths.** There are many myths about mathematics at many levels, whether they be myths about mathematical life, myths about who can and cannot do mathematics, or myths about our history and culture. A new myth is propagating throughout the mathematical community, namely the myth that the use of current AI technologies, from current companies, is inevitable. Do not believe the myth. AI tools can arise in many forms, and with many purposes. Right now, there are conversations happening within the community of research mathematicians regarding the possibility of building AI systems that are ethically produced, trained on public data of high quality, narrow in scope to ensure safety and reliability, and available for low or no cost. You do not have to be either “pro-AI” or “anti-AI” to speak the truth and reject harmful mythology. We must unlearn the myths we have been taught, so that we can view the world with clear sight. Only then can we make decisions about our future that are compassionate and just.

**Retain your humanity and sustain welcoming communities.** In response to the development of AI tools, there have been various calls to offload mathematical reasoning to the machines, shifting the role of mathematicians and the practice of doing mathematics. In this vision, how will we play the proverbial game of mathematics? Where we were previously the players, playing together on the same team, is it our future to be the coaches, with each of our teams competing against each other? Will we hoard our most promising questions? Will we retain some level of trust in each other? Will we still collaborate, and mentor, and share? Are we each destined to become our own tiny enterprise, striving to be the first to reach the next benchmark? We do not have to cede the human role in our mathematical thoughts, practices, and cultures to AI tools. We can choose a more nuanced and balanced future, one in which we emphasize and reward collaboration instead of competition, one in which we hold in the highest regard those who develop the ability to support and encourage others experiencing the authentic challenges that mathematics provides, while simultaneously holding high expectations for understanding.

**Learn about movements for the radical transformation of mathematics.** While math is often viewed as a neutral subject, as a discipline where rational thought rules the day, the reality of our history and cultures are far from this utopic vision. The reality is that our mathematical cultures and communities are the product of centuries of power struggles and political battles. Just within the United States, the founding of the Mathematical Association of America was a consequence of a closely decided vote in 1915 by a committee of the American Mathematical Society to reject the creation of a general interest mathematics journal, a journal that later became the American Mathematical Monthly. The Association for Women in Mathematics, the National Association of Mathematicians, and Spectra were founded in response to the failure of mathematicians in positions of power to recognize the dignity, value, and needs of others. Contemporary efforts such as the movement to rehumanize mathematics and the Algebra Project are important and can inform our responses in this challenging moment.

Now is not the time for rash decisions. It is not the time for decisions made with a narrow vision, applying to research alone, or teaching alone, or to any other singular focus. Now is the time for us to reach out to each other, to collaborate with intention, and to center empathy, compassion, justice, and humanity in our decisions about the future of the practices and cultures of mathematics.

*Acknowledgements*: Thanks to Drew Armstrong, Matt Beck, Laura Braun, Tendai Chitewere, Rob Davis, and Josh Grochow for their insightful comments and criticism regarding a preliminary draft of this essay.